\documentstyle[amsmath,12pt]{article}
\newtheorem{th}{Theorem}[section]

\newtheorem{prop}[th]{Proposition}
\newtheorem{cor}[th]{Corollary}
\newtheorem{defn}[th]{Definition}
\newenvironment{defn-new}{\begin{defn} \em}{\end{defn}}
\newtheorem{rem}[th]{Remark}
\newenvironment{rem-new}{\begin{rem} \em}{\end{rem}}
\newtheorem{ex}[th]{Example}
\newenvironment{ex-new}{\begin{ex} \em}{\end{ex}}

\newenvironment{notation-new}{\begin{rem} \em}{\end{rem}}

\newenvironment{agr-new}{\begin{rem} \em}{\end{rem}}

\makeatletter \@addtoreset{equation}{section} \makeatother

\makeatletter \@addtoreset{figure}{section} \makeatother

\begin{document}

\begin{center}
{\bf {\Large Comprehensive quasi-Einstein spacetime with application to
general relativity}}\\[0pt]

{\large Punam Gupta and Sanjay Kumar Singh}
\end{center}

\noindent {\bf Abstract.} The aim of this paper is to extend the notion of
all known quasi-Einstein manifolds like generalized quasi-Einstein, mixed
generalized quasi-Einstein manifold, pseudo generalized quasi-Einstein
manifold and many more and name it comprehensive quasi Einstein manifold $%
Co(QE)_{n}$. We investigate some geometric and physical properties of the
comprehensive quasi Einstein manifolds $Co(QE)_{n}$ under certain
conditions. We study the conformal and conharmonic mappings between $%
Co(QE)_{n}$ manifolds. Then we examine the $Co(QE)_{n}$ with harmonic Weyl
tensor. We define the manifold of comprehensive quasi-constant curvature and
proved that conformally flat $Co(QE)_{n}$ is manifold of comprehensive
quasi-constant curvature and vice versa. We study the general two viscous
fluid spacetime $Co(QE)_{4}$ and find out some important consequences about $%
Co(QE)_{4}$. We study $Co(QE)_{n}$ with vanishing space matter tensor.
Finally, we prove the existence of such manifolds by constructing
non-trivial example.

 \medskip

\medskip

\section{Introduction}

An Einstein manifold \cite{Besse} is a Riemannian or pseudo-Riemannian
differentiable manifold whose Ricci curvature is proportional to the metric.
Einstein manifolds named after Albert Einstein because this condition is
equivalent to saying that the metric is a solution of the vacuum Einstein
field equations with cosmological constant. Therefore they play an important
role in differential geometry as well as general theory of relativity. Four
dimensional Riemannian Einstein manifolds are also important in mathematical
physics as gravitational instantons in quantum theories of gravity.

In 2000, Chaki and Maity \cite{Chaki-Maity-00} introduced the notion of
quasi-Einstein manifolds which is a generalisation of Einstein manifold and
proved that Robertson-Walker spacetime is quasi-Einstein manifold. Also they
have some importance in the general theory of relativity.  For more details, see \cite{Chaki-Ghosal-03,De-Ghosh-04PMH,De-Ghosh-04,De-Ghosh-05,De-Ghosh-05PMD,De-De,De-Sengupta-Saha-06,Ghosh-De-Binh-06}. After that Chaki 
\cite{Chaki-01} introduced generalized quasi Einstein manifolds and find its
applications in physics. Later authors defined mixed generalized
quasi-Einstein manifold \cite{Bhatt}, nearly quasi Einstein manifold{\em \ }%
\cite{De-Gazi-08}, pseudo quasi-Einstein manifold \cite{Shaikh-09}, pseudo
generalized quasi-Einstein manifold \cite{Shaikh-Jana-08}, super
quasi-Einstein manifold \cite{Chaki-04}, mixed quasi-Einstein manifold \cite%
{Mallick-De}, mixed super quasi-Einstein manifold \cite{BTD-08},
hyper-generalized quasi-Einstein manifold \cite{Shaikh-Ozgur} and obtained
the applications of these manifolds in general theory of relativity.

In this paper, we extend the notion of generalized quasi-Einstein, mixed
generalized quasi-Einstein manifold, pseudo generalized quasi-Einstein
manifold and many more and name it comprehensive quasi Einstein manifold $%
Co(QE)_{n}$. We connect it to two fluid spacetime. Basically, when we study
the general model in theory of relativity, it needs these types of manifold
to understand the thing in better way. We study the conformal and
conharmonic mappings between $Co(QE)_{n}$ manifolds. Then we examine the $%
Co(QE)_{n}$ with harmonic Weyl tensor. We investigate geometric and physical
properties of the comprehensive quasi Einstein manifolds $Co(QE)_{n}$ under
certain conditions. We define the manifold of comprehensive quasi-constant
curvature and proved that conformally flat $Co(QE)_{n}$ is manifold of
comprehensive quasi-constant curvature and vice versa. We study the general
two viscous fluid spacetime $Co(QE)_{4}$ and find out some improtant
consequences about $Co(QE)_{4}$. We study $Co(QE)_{n}$ with vanishing space
matter tensor. Finally, we prove the existence of such manifolds by
constructing non-trivial example.

\section{Comprehensive quasi Einstein manifold}

In this section, we introduce the notion of comprehensive quasi Einstein
manifold and investigate the geometrical properties of manifold.

\begin{defn-new}
A non-flat Riemannian manifold is called a comprehensive quasi Einstein
manifold if its Ricci-tensor $S$ of type $(0,2)$ is non-zero and satisfies%
\begin{equation}
S(X,Y)=ag(X,Y)+b_{ij}\omega ^{i}(X)\omega
^{j}(Y)+c_{1}d_{1}(X,Y)+c_{2}d_{2}(X,Y),\qquad i,j=1,\ldots ,4  \label{eq-co}
\end{equation}%
for all $X,Y\in \Gamma M$, where $a,b_{ij}$,$c_{1},c_{2}$ are scalars, $%
d_{1},d_{2}$ are symmetric $(0,2)$ type tensors and $\omega ^{i}(i=1,\ldots
,4)$ are non-zero $1$-forms such that%
\[
g(X,W_{i})=\omega ^{i}(X),\qquad b_{ij}=b_{ji}
\]%
\[
g(W_{i},W_{j})=\delta _{j}^{i},\qquad {\rm trace}d_{1}=0={\rm trace}%
d_{2},\qquad d_{1}(X,W_{1})=0=d_{2}(X,W_{1}),
\]%
$W_{i}$ , $i=1,\ldots ,4$ are unit vector fields. The scalars $a,b_{ij}$,$%
c_{1},c_{2}$ are called the associated scalars, $\omega ^{i}(i=1,\ldots ,4)$
are called the associated $1$-forms and vector fields $W_{i}\left(
i=1,\ldots ,4\right) $ are called the generators of the manifold. This type
of manifold will be denoted by $Co(QE)_{n}$.
\end{defn-new}

A comprehensive quasi Einstein manifold, in particular, reduces to the
following manifolds:%
\[
\begin{tabular}{|l|l|}
\hline
Einstein \cite{Besse} & $%
\begin{array}{c}
a\not=0,\text{ all }b_{ij}=0, \\ 
c_{1}=0\text{ and }c_{2}=0%
\end{array}%
$ \\ \hline
Quasi-Einstein \cite{Chaki-Maity-00} & $%
\begin{array}{c}
a\not=0,\text{only }b_{11}\not=0 \\ 
c_{1}=0\text{ and }c_{2}=0%
\end{array}%
$ \\ \hline
generalized quasi Einstein manifold \cite{Chaki-01} & $%
\begin{array}{c}
a\not=0,\text{only }b_{11},b_{22}\not=0 \\ 
c_{1}=0\text{ and }c_{2}=0%
\end{array}%
$ \\ \hline
mixed generalized quasi-Einstein manifold \cite{Bhatt} & $%
\begin{array}{c}
a\not=0,\text{only }b_{11},b_{22},b_{12},b_{21}\not=0 \\ 
c_{1}=0\text{ and }c_{2}=0%
\end{array}%
$ \\ \hline
nearly quasi Einstein manifold{\em \ }\cite{De-Gazi-08} & $%
\begin{array}{c}
a\not=0,\text{ all }b_{ij}=0, \\ 
c_{1}\not=0\text{ and }c_{2}=0%
\end{array}%
$ \\ \hline
pseudo quasi-Einstein manifold \cite{Shaikh-09} & $%
\begin{array}{c}
a\not=0,\text{only }b_{11}\not=0 \\ 
c_{1}\not=0\text{ and }c_{2}=0%
\end{array}%
$ \\ \hline
pseudo generalized quasi-Einstein manifold \cite{Shaikh-Jana-08} & $%
\begin{array}{c}
a\not=0,\text{only }b_{11},b_{22}\not=0 \\ 
c_{1}\not=0\text{ and }c_{2}=0%
\end{array}%
$ \\ \hline
super quasi-Einstein manifold \cite{Chaki-04} & $%
\begin{array}{c}
a\not=0,\text{only }b_{11},b_{12},b_{21}\not=0 \\ 
c_{1}\not=0\text{ and }c_{2}=0%
\end{array}%
$ \\ \hline
mixed quasi-Einstein manifold \cite{Mallick-De} & $%
\begin{array}{c}
a\not=0,\text{only }b_{12},b_{21}\not=0 \\ 
c_{1}=0\text{ and }c_{2}=0%
\end{array}%
$ \\ \hline
mixed super quasi-Einstein manifold \cite{BTD-08} & $%
\begin{array}{c}
a\not=0,\text{only }b_{11},b_{22},b_{12},b_{21}\not=0 \\ 
c_{1}\not=0\text{ and }c_{2}=0%
\end{array}%
$ \\ \hline
hyper-generalized quasi-Einstein manifold \cite{Shaikh-Ozgur} & $%
\begin{array}{c}
a\not=0,\text{only }b_{11},b_{12},b_{21},b_{13},b_{31}\not=0 \\ 
c_{1}=0\text{ and }c_{2}=0%
\end{array}%
$ \\ \hline
\end{tabular}%
\ 
\]%
Beacuse all the known quasi-Einstein manifolds are the particular case of $%
Co(QE)_{n}$, justifies the name of this manifold.

Let $\{e_{i}:i=1,...,n\}$ be an orthonormal basis of the tangent space at
any point of the manifold. Then setting $X=Y=e_{i}$ in (\ref{eq-co}) and
taking summation over $i$, $1\leq i\leq n,$ we obtain 
\begin{equation}
r=an+b_{11}+b_{22}+b_{33}+b_{44}.  \label{eq-co-r}
\end{equation}
It is easy to find that 
\[
S(W_{1},W_{1})=a+b_{11}, 
\]%
\[
S(W_{2},W_{2})=a+b_{22}+c_{1}d_{1}(W_{2},W_{2})+c_{2}d_{2}(W_{2},W_{2}), 
\]%
\begin{eqnarray*}
S(W_{3},W_{3}) &=&a+b_{33}+c_{1}d_{1}(W_{3},W_{3})+c_{2}d_{2}(W_{3},W_{3}),
\\
S(W_{4},W_{4}) &=&a+b_{44}+c_{1}d_{1}(W_{4},W_{4})+c_{2}d_{2}(W_{4},W_{4}),
\\
S(W_{1},W_{2}) &=&b_{12}, \\
S(W_{2},W_{3}) &=&b_{23}+c_{1}d_{1}(W_{2},W_{3})+c_{2}d_{2}(W_{2},W_{3}), \\
S(W_{3},W_{4}) &=&b_{34}+c_{1}d_{1}(W_{3},W_{4})+c_{2}d_{2}(W_{3},W_{4}), \\
S(W_{3},W_{1}) &=&b_{31}, \\
S(W_{4},W_{1}) &=&b_{41}, \\
S(W_{2},W_{4}) &=&b_{24}+c_{1}d_{1}(W_{2},W_{4})+c_{2}d_{2}(W_{2},W_{4}).
\end{eqnarray*}%
Recall that $S(W,W)$ is the Ricci curvature in the direction of $W$ if $W$
is a unit vector field. Therefore, we can state that

\begin{th}
In a $Co(QE)_{n}(n>2)$, the scalars $a+b_{11}$, $%
a+b_{22}+c_{1}d_{1}(W_{2},W_{2})+c_{2}d_{2}(W_{2},W_{2})$, $%
a+b_{33}+c_{1}d_{1}(W_{3},W_{3})+c_{2}d_{2}(W_{3},W_{3})$ and $%
a+b_{44}+c_{1}d_{1}(W_{4},W_{4})+c_{2}d_{2}(W_{4},W_{4})$ are the Ricci
curvatures in the directions of the generators $W_{1},W_{2}$, $W_{3}$ and $%
W_{4}$, respectively.
\end{th}

Let $s^{2}$and $t_{1}^{2},t_{2}^{2}$ denote the sqaures of the length of the
Ricci tensor $S$ and the structure tensors $d_{1},d_{2}$, respectively, that
is, 
\[
s^{2}=\sum\limits_{i=1}^{n}S(Qe_{i},e_{i}), 
\]%
\[
t_{1}^{2}=\sum\limits_{i=1}^{n}d_{1}(D_{1}e_{i},e_{i}), 
\]%
\[
t_{2}^{2}=\sum\limits_{i=1}^{n}d_{2}(D_{2}e_{i},e_{i}), 
\]%
where $Q$, $D_{1}$, $D_{2}$ are symmetric endomorphism of the tangent space
at each point corresponding to the Ricci tensor $S$, $d_{1}$, $d_{2}$,
respectively. Then it is easy to compute that 
\begin{eqnarray*}
s^{2}-c_{1}t_{1}^{2}-c_{2}t_{2}^{2}
&=&na^{2}+b_{11}^{2}+b_{22}^{2}+b_{33}^{2}+b_{44}^{2}+2a\left(
b_{11}+b_{22}+b_{33}+b_{44}\right) \\
&&+2b_{12}^{2}+2b_{23}^{2}+2b_{14}^{2}+2b_{34}^{2}+2b_{13}^{2}+2b_{24}^{2} \\
&&+2c_{1}\left( 
\begin{array}{c}
b_{22}d_{1}(W_{2},W_{2})+b_{33}d_{1}(W_{3},W_{3})+b_{44}d_{1}(W_{4},W_{4})
\\ 
+2b_{23}d_{1}(W_{2},W_{3})+2b_{34}d_{1}(W_{4},W_{3})+2b_{24}d_{1}(W_{2},W_{4})%
\end{array}%
\right) \\
&&+2c_{2}\left( 
\begin{array}{c}
b_{22}d_{2}(W_{2},W_{2})+b_{33}d_{2}(W_{3},W_{3})+b_{44}d_{2}(W_{4},W_{4})
\\ 
+2b_{23}d_{2}(W_{2},W_{3})+2b_{34}d_{2}(W_{4},W_{3})+2b_{24}d_{2}(W_{2},W_{4})%
\end{array}%
\right) \\
&&+\left( c_{1}+c_{2}\right) \sum\limits_{i=1}^{n}g(D_{1}e_{i},D_{2}e_{i}).
\end{eqnarray*}%
Now, consider functions $a,b_{ij}$,$c_{1},c_{2}$ are constant. Then 
\begin{eqnarray}
(\nabla _{W}S)(X,Y) &=&b_{11}((\nabla _{W}\omega ^{1})(X)\omega
^{1}(Y)+\omega ^{1}(X)(\nabla _{W}\omega ^{1})(Y))  \nonumber \\
&&+b_{22}((\nabla _{W}\omega ^{2})(X)\omega ^{2}(Y)+\omega ^{2}(X)(\nabla
_{W}\omega ^{2})(Y))  \nonumber \\
&&+b_{33}((\nabla _{W}\omega ^{3})(X)\omega ^{3}(Y)+\omega ^{3}(X)(\nabla
_{W}\omega ^{3})(Y))  \nonumber \\
&&+b_{44}((\nabla _{W}\omega ^{4})(X)\omega ^{4}(Y)+\omega ^{4}(X)(\nabla
_{W}\omega ^{4})(Y))  \nonumber \\
&&+b_{12}((\nabla _{W}\omega ^{1})(X)\omega ^{2}(Y)+\omega ^{1}(X)(\nabla
_{W}\omega ^{2})(Y))  \nonumber \\
&&+b_{12}((\nabla _{W}\omega ^{1})(Y)\omega ^{2}(X)+\omega ^{1}(Y)(\nabla
_{W}\omega ^{2})(X))  \nonumber \\
&&+b_{13}((\nabla _{W}\omega ^{1})(X)\omega ^{3}(Y)+\omega ^{1}(X)(\nabla
_{W}\omega ^{3})(Y))  \nonumber \\
&&+b_{13}((\nabla _{W}\omega ^{1})(Y)\omega ^{3}(X)+\omega ^{1}(Y)(\nabla
_{W}\omega ^{3})(X))  \nonumber \\
&&+b_{23}((\nabla _{W}\omega ^{2})(X)\omega ^{3}(Y)+\omega ^{2}(X)(\nabla
_{W}\omega ^{3})(Y))  \nonumber \\
&&+b_{23}((\nabla _{W}\omega ^{2})(Y)\omega ^{3}(X)+\omega ^{2}(Y)(\nabla
_{W}\omega ^{3})(X))  \nonumber \\
&&+b_{24}((\nabla _{W}\omega ^{2})(X)\omega ^{4}(Y)+\omega ^{2}(X)(\nabla
_{W}\omega ^{4})(Y))  \nonumber \\
&&+b_{24}((\nabla _{W}\omega ^{2})(Y)\omega ^{4}(X)+\omega ^{2}(Y)(\nabla
_{W}\omega ^{4})(X))  \nonumber \\
&&+b_{34}((\nabla _{W}\omega ^{3})(X)\omega ^{4}(Y)+\omega ^{3}(X)(\nabla
_{W}\omega ^{4})(Y))  \nonumber \\
&&+b_{43}((\nabla _{W}\omega ^{3})(Y)\omega ^{4}(X)+\omega ^{3}(Y)(\nabla
_{W}\omega ^{4})(X))  \nonumber \\
&&+b_{14}((\nabla _{W}\omega ^{4})(X)\omega ^{1}(Y)+\omega ^{4}(X)(\nabla
_{W}\omega ^{1})(Y))  \nonumber \\
&&+b_{14}((\nabla _{W}\omega ^{4})(Y)\omega ^{1}(Y)+\omega ^{4}(Y)(\nabla
_{W}\omega ^{1})(X))  \nonumber \\
&&+c_{1}(\nabla _{W}d_{1})(X,Y)+c_{2}(\nabla _{W}d_{2})(X,Y). \label{eq-S}
\end{eqnarray}%
Next, we give some definitions for further use:

\begin{defn-new}
A non-flat $n$-dimensional Riemannian manifold $(M,g),(n>3)$ is called a
semi-pseudo Ricci symmetric manifold \cite{Tarafdar} if the Ricci tensor $S$
of type $(0,2)$ is non-zero and satisfies the condition 
\begin{equation}
\left( \nabla _{X}S\right) (Y,Z)=\pi (Y)S(X,Z)+\pi (Z)S(X,Y),  \label{eq-sp}
\end{equation}%
where $\nabla $ denotes the Levi-Civita connection and $\pi $ is a non-zero $%
1$-form such that $g(X,\Pi )=\pi (X)$ for all vector fields $X$, $\Pi $
being the vector field corresponding to the associated $1$-form $\pi $. If $%
\pi =0$, then the manifold is called Ricci symmetric.
\end{defn-new}

\begin{defn-new}
Let $(M,g)$ be an $n$-dimensional Riemannian manifold. The Ricci tensor $S$
is called cyclic parallel \cite{Gray} if it satisfies the condition%
\begin{equation}
(\nabla _{X}S)(Y,Z)+(\nabla _{Y}S)(Z,X)+(\nabla _{Z}S)(X,Y)=0,  \label{eq-cp}
\end{equation}%
for any vector fields $X,Y,Z$ on $M$. Since every Einstein manifold
satisfies (\ref{eq-cp}), therefore the manifold with  parallel Ricci
tensor is also known as Einstein-like manifolds.
\end{defn-new}

\begin{defn-new}
A symmetric tensor field $A$ of type $(0,2)$ on a Riemannian manifold $(M,g)$
is said to be a Codazzi tensor \cite{Besse} if it satisfies the condition 
\[
(\nabla _{X}A)(Y,Z)=(\nabla _{Y}A)(X,Z), 
\]%
for any vector fields $X,Y,Z$ on $M$.
\end{defn-new}

\begin{defn-new}
A $\varphi $(Ric)-vector field \cite{Hen,Kirik} is a vector field on an $n$%
-dimensional Riemannian manifold $M$ and Levi-Civita
connection $\nabla $, which satisfies the condition 
\[
\nabla \varphi =\mu S, 
\]%
where $\mu $ is a constant and $S$ is the Ricci tensor. If $M$ is an
Einstein manifold, the vector field $\varphi $ is concircular. If $\mu
\not=0 $, then the vector field $\varphi $ is proper $\varphi $(Ric)-vector
field. When $\mu =0$, the vector field $\varphi $ is covariantly constant.
\end{defn-new}

In other words, we define

\begin{defn-new}
A vector field $U$ is said to be concircular \cite{Schouten} if
\[
\nabla _{X}U=\rho X, 
\]%
where $\rho $ is a function on the manifold. If $\rho $ is a non-zero
constant, then vector field $U$ is said to be concurrent \cite{Schouten}. If 
$\rho =0$, the vector field reduces to a parallel vector field.
\end{defn-new}

\begin{defn-new}
A vector field $W$ corresponding to the associated $1$-form $\omega $ is
said to be recurrent if \cite{Schouten}%
\begin{equation}
(\nabla _{X}\omega )(Y)=\phi (X)\omega (Y),  \label{eq-rec}
\end{equation}%
where $\phi $ is a non-zero $1$-form.
\end{defn-new}

\begin{defn-new}
A non- flat Riemannian manifold is said to be generalized Ricci recurrent 
\cite{De-Guha} if its Ricci tensor $S$ of type $(0,2)$ satisfies the
condition 
\[
(\nabla _{X}S)(Y,Z)=\alpha (X)S(Y,Z)+\beta (X)S(Y,Z), 
\]%
where $\alpha (X)$,$\beta (X)$ are non-zero $1$-forms. If $\beta (X)=0$,
then it reduces to Ricci recurrent manifold \cite{Patterson-52}.
\end{defn-new}

\section{Conformal \& Conharmonic mappings of $Co(QE)_{n}$}

In this section, we consider the conformal and conharmonic mappings between $%
Co(QE)_{n}$. Let $M$ and $N$ be two $Co(QE)_{n}$ with metrics $g$ and $%
\tilde{g}$, respectively. For this, we need the following result which
directly comes by use the result of \cite[Th 2.3]{Kirik}:

\begin{th}
Consider $Co(QE)_{n}$ $(n>3)$ such that the associated scalars are
constants. If $Co(QE)_{n}$ admits a $\varphi $(Ric)-vector field, then the
length of $\varphi $ is constant.
\end{th}

Now, we state some definitions:

\begin{defn-new}
A diffeomorphism $f:(M,g)\rightarrow (N,\tilde{g})$ is said to be conformal
mapping \cite{Eisenhart} if 
\[
\tilde{g}=e^{2\sigma }g, 
\]
where $\sigma $ is a function on $M$. If $\sigma $ is constant, then
conformal mapping is called homothetic mapping.
\end{defn-new}

\begin{defn-new}
A conformal mapping $f:(M,g)\rightarrow (N,\tilde{g})$ is said to be
conharmonic mapping (transformation) \cite{Ishii} if it satisfies 
\[
\triangle \sigma =-\frac{n-2}{2}\left\Vert {grad}\sigma \right\Vert ^{2}, 
\]%
where $\triangle $ is Laplace-Beltrami operator.
\end{defn-new}

Now, we are using the result of \cite[Th 4.3]{Kirik}, which is also true for 
$Co(QE)_{n}$. So we can state that

\begin{th}
Let $f:M\rightarrow N$ be conformal mapping between $Co(QE)_{n}$ manifolds
such that Ricci tensors of $M$ and $N$ are Codazzi type. If the vector field
generated by the $1$-form $\sigma $ is a $\sigma $(Ric)-vector field, then
either this conformal mapping is homothetic or satisfies 
\[
\mu =\frac{(2-n)(n-1)\left\Vert {grad}\sigma \right\Vert ^{2}-r}{2(n-1)r}, 
\]
$\mu $ denotes the constant corresponding to the $\sigma $(Ric)-vector field
and $r\not=0$.
\end{th}

\begin{th}
Let $f:M\rightarrow N$ be conformal mapping between $Co(QE)_{n}$ manifolds
such that Ricci tensors of $M$ and $N$ are Codazzi type. If ${grad}\sigma $
is a concircular vector field, then either $W_{i}$ $(i=1,2,3,4)$ and ${grad}%
\sigma $ are orthogonal or $b_{12}=b_{13}=b_{14}=0$ and $\rho =\dfrac{%
(n-2)(1-n)\triangle _{1}\sigma -b_{22}-b_{33}-b_{44}}{(n+2)(n-1)}$, where $%
\rho $ and $\triangle _{1}\sigma $ denote the function corresponding to the
concircular vector field and first Beltrami's symbol, respectively.
\end{th}

\noindent {\bf Proof.} The proof is similar to the proof of \cite[Th 4.4%
]{Kirik}.

\begin{th}
Let $f:M\rightarrow N$ be conformal mapping between $Co(QE)_{n}$ manifolds.
Then conformal mapping is conharmonic if and only if the associated scalars $%
\tilde{a}$,$\tilde{b}_{11}$,$\tilde{b}_{22},\tilde{b}_{33},\tilde{b}_{44}$
be transformed by $\tilde{a}=e^{-2\sigma }a,$ $\tilde{b}_{11}=e^{-2\sigma
}b_{11}$, $\tilde{b}_{22}=e^{-2\sigma }b_{22}$, $\tilde{b}_{33}=e^{-2\sigma
}b_{33}$, $\tilde{b}_{44}$=$e^{-2\sigma }b_{44}$.
\end{th}

\noindent {\bf Proof.} The proof is similar to the proof of \cite[Th 4.5%
]{Kirik}.

\section{$Co(QE)_{n}(n>3)$ with harmonic Weyl tensor}

The Weyl tensor (conformal curvature tensor) \cite{Neill} is invariant under
conformal mapping and is given as 
\begin{eqnarray}
{\cal C}(X,Y,Z,W) &=&R(X,Y,Z,W)-\frac{1}{n-2}\left(
S(X,W)g(Y,Z)-S(Y,W)g(X,Z)\right.  \nonumber \\
&&+\left. g(X,W)S(Y,Z)-g(Y,W)S(X,Z)\right)  \nonumber \\
&&+\frac{r}{(n-1)(n-2)}\left( g(X,W)g(Y,Z)-g(Y,W)g(X,Z)\right) .
\label{eq-cct}
\end{eqnarray}%
Using above equation, we have 
\begin{eqnarray}
({\rm div}{\cal C})(X,Y,Z) &=&\frac{n-3}{n-2}\left( (\nabla
_{X}S)(Y,Z)-(\nabla _{Y}S)(X,Z)\right)  \nonumber \\
&&+\frac{n-3}{(n-1)\left( n-2\right) }\left( (\nabla _{Y}r)g(X,Z)-(\nabla
_{X}r)g(Y,Z)\right).\label{eq-div-1}
\end{eqnarray}%
Weyl tensor is said to be harmonic if the divergence of ${\cal C}$ vanishes.
In $3$-dim, this condition is equivalent to local conformally flatness.
Nevertheless, when $n>3$, harmonic Weyl tensor is a weaker condition since
locally conformally flatness is equivalent to the vanishing of the Weyl
tensor.

By using (\ref{eq-cct}), we get 
\[
({\rm div}{\cal C})(X,Y,Z)=-\frac{(n-3)}{(n-2)}C(X,Y,Z), 
\]%
where $C(X,Y,Z)$ is the cotton tensor given by 
\begin{equation}
C(X,Y,Z)=(\nabla _{Z}S)(X,Y)-(\nabla _{Y}S)(Z,X)-\frac{1}{2(n-1)}\left(
(\nabla _{Z}r)g(Y,Z)-(\nabla _{Y}r)g(X,Z)\right) .  \label{eq-cot}
\end{equation}%
If $n>4$, harmonic Weyl tensor is equivalent to the vanishing of the Cotton
tensor.

Consider a $Co(QE)_{n}(n>3)$. If all the scalars $a,b_{ij},c_{1},c_{2}$ are
constant, then the covariant derivative of $r$ becomes zero. Now if we
consider that the generators $W_{i}$ of the manifold are recurrent vector
field with associated $1$-forms $\omega ^{i}$, respectively, not being the $%
1 $-form of recurrence, gives $\nabla _{X}W_{i}=\pi _{i}(X)W_{i}$, where $%
\pi _{i}$ are the $1$-form of recurrence, we get%
\begin{equation}
g(\nabla _{X}W_{i},Y)=\pi _{i}(X)g(W_{i},Y),\qquad (\nabla _{X}\omega
^{i})(Y)=\pi _{i}(X)\omega ^{i}(Y).  \label{eq-1}
\end{equation}%
By (\ref{eq-1}), we have $(\nabla _{X}\omega ^{i})(W_{i})=$ $g(\nabla
_{X}W_{i},W_{i})=0$. Therefore $\pi _{i}(X)=0$ for all $X$. Also, assume
that the structure tensors $d_{1},d_{2}$ are of Codazzi type, then using (%
\ref{eq-S}) in (\ref{eq-div-1}), we get $({\rm div}{\cal C})(X,Y,Z)=0$.

Thus we can state the following:

\begin{th}
If in a $Co(QE)_{n}$ $(n>3),$ the associated scalars are constants and
generators $W_{i}$ of the manifold are recurrent vector fields with the
associated $1$-form $\omega ^{i},$ respectively, not being the the $1$-form
of recurrence and the structure tensors $d_{1},d_{2}$ are of Codazzi type,
then the Weyl tensor of the manifold is harmonic.
\end{th}

\section{Conformally flat $Co(QE)_{n}(n>3)$}

Let $M$ be a conformally flat $Co(QE)_{n}(n>3)$. By using (\ref{eq-cct}),
curvature tensor $R$ is 
\begin{eqnarray}
R(X,Y,Z,W) &=&\frac{1}{n-2}\left( S(X,W)g(Y,Z)-S(Y,W)g(X,Z)\right.  \nonumber
\\
&&\left. +g(X,W)S(Y,Z)-g(Y,W)S(X,Z)\right)  \nonumber \\
&&-\frac{r}{(n-1)(n-2)}\left( g(X,W)g(Y,Z)-g(Y,W)g(X,Z)\right) .
\label{eq-con}
\end{eqnarray}%
By using (\ref{eq-co}) and (\ref{eq-co-r}) in (\ref{eq-con}), we get 
\begin{eqnarray*}
R(X,Y,Z,W) &=&-\frac{(a+b_{11}+b_{22}+b_{33}+b_{44})}{(n-1)(n-2)}\left(
g(X,W)g(Y,Z)-g(Y,W)g(X,Z)\right) \\
&&+\frac{c_{1}}{n-2}\left( d_{1}(X,W)g(Y,Z)-d_{1}(Y,W)g(X,Z)\right. \\
&&\left. +d_{1}(Y,Z)g(X,W)-d_{1}(X,Z)g(Y,W)\right) \\
&&+\frac{c_{2}}{n-2}\left( d_{2}(X,W)g(Y,Z)-d_{2}(Y,W)g(X,Z)\right. \\
&&\left. +d_{2}(Y,Z)g(X,W)-d_{2}(X,Z)g(Y,W)\right) \\
&&+\frac{b_{11}}{n-2}\left( \omega ^{1}(X)\omega ^{1}(W)g(Y,Z)-\omega
^{1}(Y)\omega ^{1}(W)g(X,Z)\right. \\
&&\left. +\omega ^{1}(Y)\omega ^{1}(Z)g(X,W)-\omega ^{1}(X)\omega
^{1}(Z)g(Y,W)\right) \\
&&+\frac{b_{22}}{n-2}\left( \omega ^{2}(X)\omega ^{2}(W)g(Y,Z)-\omega
^{2}(Y)\omega ^{2}(W)g(X,Z)\right. \\
&&\left. +\omega ^{2}(Y)\omega ^{2}(Z)g(X,W)-\omega ^{2}(X)\omega
^{2}(Z)g(Y,W)\right) \\
&&+\frac{b_{33}}{n-2}\left( \omega ^{3}(X)\omega ^{3}(W)g(Y,Z)-\omega
^{3}(Y)\omega ^{3}(W)g(X,Z)\right. \\
&&\left. +\omega ^{3}(Y)\omega ^{3}(Z)g(X,W)-\omega ^{3}(X)\omega
^{3}(Z)g(Y,W)\right) \\
&&+\frac{b_{44}}{n-2}\left( \omega ^{4}(X)\omega ^{4}(W)g(Y,Z)-\omega
^{4}(Y)\omega ^{4}(W)g(X,Z)\right. \\
&&\left. +\omega ^{4}(Y)\omega ^{4}(Z)g(X,W)-\omega ^{4}(X)\omega
^{4}(Z)g(Y,W)\right) \\
&&+\frac{b_{12}}{n-2}\left( 
\begin{array}{c}
(\omega ^{1}(X)\omega ^{2}(W)+\omega ^{2}(X)\omega ^{1}(W))g(Y,Z) \\ 
-(\omega ^{1}(Y)\omega ^{2}(W)+\omega ^{2}(Y)\omega ^{1}(W))g(X,Z) \\ 
+(\omega ^{1}(Y)\omega ^{2}(Z)+\omega ^{2}(Y)\omega ^{1}(Z))g(X,W) \\ 
-(\omega ^{1}(X)\omega ^{2}(Z)+\omega ^{2}(X)\omega ^{1}(Z))g(Y,W)%
\end{array}%
\right) \\
&&+\frac{b_{14}}{n-2}\left( 
\begin{array}{c}
(\omega ^{1}(X)\omega ^{4}(W)+\omega ^{4}(X)\omega ^{1}(W))g(Y,Z) \\ 
-(\omega ^{1}(Y)\omega ^{4}(W)+\omega ^{4}(Y)\omega ^{1}(W))g(X,Z) \\ 
+(\omega ^{1}(Y)\omega ^{4}(Z)+\omega ^{4}(Y)\omega ^{1}(Z))g(X,W) \\ 
-(\omega ^{1}(X)\omega ^{4}(Z)+\omega ^{4}(X)\omega ^{1}(Z))g(Y,W)%
\end{array}%
\right) \\
&&+\frac{b_{23}}{n-2}\left( 
\begin{array}{c}
(\omega ^{2}(X)\omega ^{3}(W)+\omega ^{3}(X)\omega ^{2}(W))g(Y,Z) \\ 
-(\omega ^{2}(Y)\omega ^{3}(W)+\omega ^{3}(Y)\omega ^{2}(W))g(X,Z) \\ 
+(\omega ^{2}(Y)\omega ^{3}(Z)+\omega ^{3}(Y)\omega ^{2}(Z))g(X,W) \\ 
-(\omega ^{2}(X)\omega ^{3}(Z)+\omega ^{3}(X)\omega ^{2}(Z))g(Y,W)%
\end{array}%
\right) \\
&&+\frac{b_{24}}{n-2}\left( 
\begin{array}{c}
(\omega ^{2}(X)\omega ^{4}(W)+\omega ^{4}(X)\omega ^{2}(W))g(Y,Z) \\ 
-(\omega ^{2}(Y)\omega ^{4}(W)+\omega ^{4}(Y)\omega ^{2}(W))g(X,Z) \\ 
+(\omega ^{2}(Y)\omega ^{4}(Z)+\omega ^{4}(Y)\omega ^{2}(Z))g(X,W) \\ 
-(\omega ^{2}(X)\omega ^{4}(Z)+\omega ^{4}(X)\omega ^{2}(Z))g(Y,W)%
\end{array}%
\right) \\
&&+\frac{b_{31}}{n-2}\left( 
\begin{array}{c}
(\omega ^{3}(X)\omega ^{1}(W)+\omega ^{1}(X)\omega ^{3}(W))g(Y,Z) \\ 
-(\omega ^{3}(Y)\omega ^{1}(W)+\omega ^{1}(Y)\omega ^{3}(W))g(X,Z) \\ 
+(\omega ^{3}(Y)\omega ^{1}(Z)+\omega ^{1}(Y)\omega ^{3}(Z))g(X,W) \\ 
-(\omega ^{3}(X)\omega ^{1}(Z)+\omega ^{1}(X)\omega ^{3}(Z))g(Y,W)%
\end{array}%
\right) \\
&&+\frac{b_{34}}{n-2}\left( 
\begin{array}{c}
(\omega ^{3}(X)\omega ^{4}(W)+\omega ^{4}(X)\omega ^{3}(W))g(Y,Z) \\ 
-(\omega ^{3}(Y)\omega ^{4}(W)+\omega ^{4}(Y)\omega ^{3}(W))g(X,Z) \\ 
+(\omega ^{3}(Y)\omega ^{4}(Z)+\omega ^{4}(Y)\omega ^{3}(Z))g(X,W) \\ 
-(\omega ^{3}(X)\omega ^{4}(Z)+\omega ^{4}(X)\omega ^{3}(Z))g(Y,W)%
\end{array}%
\right).
\end{eqnarray*}
Now, we can define

\begin{defn-new}
A Riemannian manifold $M(n>2)$ is said to be of comprehensive quasi-constant
curvature if it is conformally flat and its satisfies 
\begin{eqnarray}
R(X,Y,Z,W) &=&a_{1}\left( g(X,W)g(Y,Z)-g(Y,W)g(X,Z)\right)  \nonumber \\
&&+a_{2}\left( 
\begin{array}{c}
d_{1}(X,W)g(Y,Z)-d_{1}(Y,W)g(X,Z) \\ 
+d_{1}(Y,Z)g(X,W)-d_{1}(X,Z)g(Y,W)%
\end{array}%
\right)  \nonumber \\
&&+a_{3}\left( 
\begin{array}{c}
d_{2}(X,W)g(Y,Z)-d_{2}(Y,W)g(X,Z) \\ 
+d_{2}(Y,Z)g(X,W)-d_{2}(X,Z)g(Y,W)%
\end{array}%
\right)  \nonumber \\
&&+a_{4}\left( 
\begin{array}{c}
\omega ^{1}(X)\omega ^{1}(W)g(Y,Z)-\omega ^{1}(Y)\omega ^{1}(W)g(X,Z) \\ 
+\omega ^{1}(Y)\omega ^{1}(Z)g(X,W)-\omega ^{1}(X)\omega ^{1}(Z)g(Y,W)%
\end{array}%
\right)  \nonumber \\
&&+a_{5}\left( 
\begin{array}{c}
\omega ^{2}(X)\omega ^{2}(W)g(Y,Z)-\omega ^{2}(Y)\omega ^{2}(W)g(X,Z) \\ 
+\omega ^{2}(Y)\omega ^{2}(Z)g(X,W)-\omega ^{2}(X)\omega ^{2}(Z)g(Y,W)%
\end{array}%
\right)  \nonumber \\
&&+a_{6}\left( 
\begin{array}{c}
\omega ^{3}(X)\omega ^{3}(W)g(Y,Z)-\omega ^{3}(Y)\omega ^{3}(W)g(X,Z) \\ 
+\omega ^{3}(Y)\omega ^{3}(Z)g(X,W)-\omega ^{3}(X)\omega ^{3}(Z)g(Y,W)%
\end{array}%
\right)  \nonumber \\
&&+a_{7}\left( 
\begin{array}{c}
\omega ^{4}(X)\omega ^{4}(W)g(Y,Z)-\omega ^{4}(Y)\omega ^{4}(W)g(X,Z) \\ 
+\omega ^{4}(Y)\omega ^{4}(Z)g(X,W)-\omega ^{4}(X)\omega ^{4}(Z)g(Y,W)%
\end{array}%
\right)  \nonumber \\
&&+a_{8}\left( 
\begin{array}{c}
(\omega ^{1}(X)\omega ^{2}(W)+\omega ^{2}(X)\omega ^{1}(W))g(Y,Z) \\ 
-(\omega ^{1}(Y)\omega ^{2}(W)+\omega ^{2}(Y)\omega ^{1}(W))g(X,Z) \\ 
+(\omega ^{1}(Y)\omega ^{2}(Z)+\omega ^{2}(Y)\omega ^{1}(Z))g(X,W) \\ 
-(\omega ^{1}(X)\omega ^{2}(Z)+\omega ^{2}(X)\omega ^{1}(Z))g(Y,W)%
\end{array}%
\right)  \nonumber \\
&&+a_{9}\left( 
\begin{array}{c}
(\omega ^{1}(X)\omega ^{4}(W)+\omega ^{4}(X)\omega ^{1}(W))g(Y,Z) \\ 
-(\omega ^{1}(Y)\omega ^{4}(W)+\omega ^{4}(Y)\omega ^{1}(W))g(X,Z) \\ 
+(\omega ^{1}(Y)\omega ^{4}(Z)+\omega ^{4}(Y)\omega ^{1}(Z))g(X,W) \\ 
-(\omega ^{1}(X)\omega ^{4}(Z)+\omega ^{4}(X)\omega ^{1}(Z))g(Y,W)%
\end{array}%
\right)  \nonumber \\
&&+a_{10}\left( 
\begin{array}{c}
(\omega ^{2}(X)\omega ^{3}(W)+\omega ^{3}(X)\omega ^{2}(W))g(Y,Z) \\ 
-(\omega ^{2}(Y)\omega ^{3}(W)+\omega ^{3}(Y)\omega ^{2}(W))g(X,Z) \\ 
+(\omega ^{2}(Y)\omega ^{3}(Z)+\omega ^{3}(Y)\omega ^{2}(Z))g(X,W) \\ 
-(\omega ^{2}(X)\omega ^{3}(Z)+\omega ^{3}(X)\omega ^{2}(Z))g(Y,W)%
\end{array}%
\right)  \nonumber \\
&&+a_{11}\left( 
\begin{array}{c}
(\omega ^{2}(X)\omega ^{4}(W)+\omega ^{4}(X)\omega ^{2}(W))g(Y,Z) \\ 
-(\omega ^{2}(Y)\omega ^{4}(W)+\omega ^{4}(Y)\omega ^{2}(W))g(X,Z) \\ 
+(\omega ^{2}(Y)\omega ^{4}(Z)+\omega ^{4}(Y)\omega ^{2}(Z))g(X,W) \\ 
-(\omega ^{2}(X)\omega ^{4}(Z)+\omega ^{4}(X)\omega ^{2}(Z))g(Y,W)%
\end{array}%
\right)  \nonumber \\
&&+a_{12}\left( 
\begin{array}{c}
(\omega ^{3}(X)\omega ^{1}(W)+\omega ^{1}(X)\omega ^{3}(W))g(Y,Z) \\ 
-(\omega ^{3}(Y)\omega ^{1}(W)+\omega ^{1}(Y)\omega ^{3}(W))g(X,Z) \\ 
+(\omega ^{3}(Y)\omega ^{1}(Z)+\omega ^{1}(Y)\omega ^{3}(Z))g(X,W) \\ 
-(\omega ^{3}(X)\omega ^{1}(Z)+\omega ^{1}(X)\omega ^{3}(Z))g(Y,W)%
\end{array}%
\right)  \nonumber \\
&&+a_{13}\left( 
\begin{array}{c}
(\omega ^{3}(X)\omega ^{4}(W)+\omega ^{4}(X)\omega ^{3}(W))g(Y,Z) \\ 
-(\omega ^{3}(Y)\omega ^{4}(W)+\omega ^{4}(Y)\omega ^{3}(W))g(X,Z) \\ 
+(\omega ^{3}(Y)\omega ^{4}(Z)+\omega ^{4}(Y)\omega ^{3}(Z))g(X,W) \\ 
-(\omega ^{3}(X)\omega ^{4}(Z)+\omega ^{4}(X)\omega ^{3}(Z))g(Y,W)%
\end{array}%
\right),  \label{eq-com-cur}
\end{eqnarray}%
where $a_{1},\ldots a_{13}$ are non-zero scalars, $\omega ^{i}(i=1,2,3,4)$
are $1$-forms and $d_{1},d_{2}$ are symmetric $(0,2)$ type trace free
tensor.
\end{defn-new}

Then we have the following particular cases:

\begin{enumerate}
\item It is said to be of constant curvature \cite{Yano} if $a_{1}\not=0$
and $a_{2}=\cdots =a_{13}=0.$

\item It is said to be of quasi-constant curvature \cite{Chen-Yano} if $%
a_{1}\not=0\not=a_{4}$ and $a_{2}=a_{3}=a_{5}=\cdots =a_{13}=0.$

\item It is said to be of generalized quasi-constant curvature \cite%
{De-Ghosh-04.} if $a_{1}\not=0,a_{4}\not=0,a_{5}\not=0$ and $%
a_{2}=a_{3}=a_{6}=\cdots =a_{13}=0.$

\item It is said to be of pseudo quasi-constant curvature \cite{Shaikh-09}
if $a_{1}\not=0,a_{2}\not=0,a_{4}\not=0$ and $a_{3}=a_{5}=\cdots =a_{13}=0.$

\item It is said to be of pseudo generalized quasi-constant curvature \cite%
{Shaikh-Jana-08} if $a_{1}\not=0,a_{2}\not=0,a_{4}\not=0,a_{5}\not=0$ and $%
a_{3}=a_{6}=\cdots =a_{13}=0.$

\item It is said to be of mixed quasi-constant curvature \cite{Mallick-17}
if $a_{1}\not=0,a_{8}\not=0$, $a_{2}=\cdots =a_{7}=a_{9}=\cdots =a_{13}=0$.

\item It is said to be of super quasi-constant curvature \cite{Chaki-04} if $%
a_{1}\not=0,a_{2}\not=0,a_{4}\not=0,a_{8}\not=0$ and $%
a_{3}=a_{5}=a_{6}=a_{7}=a_{9}=\cdots =a_{13}=0.$

\item It is said to be of mixed super quasi-constant curvature \cite{BTD-08}
if $a_{1}\not=0,a_{2}\not=0,a_{4}\not=0,a_{5}\not=0,a_{8}\not=0,$ and $%
a_{3}=a_{6}=a_{7}=a_{9}=\cdots =a_{13}=0.$

\item It is said to be of nearly quasi-constant curvature \cite{De-Gazi-08}
if $a_{1}\not=0,a_{2}\not=0$ and $a_{3}=\cdots =a_{13}=0.$

\item It is said to be of mixed generalized quasi-constant curvature \cite%
{Bhatt} if $a_{1}\not=0,a_{4}\not=0,a_{5}\not=0,a_{8}\not=0,$ and $%
a_{2}=a_{3}=a_{6}=a_{7}=a_{9}=\cdots =a_{13}=0.$

\item It is said to be of hyper-generalized quasi-constant curvature \cite%
{Shaikh-Ozgur} if $a_{1}\not=0,a_{4}\not=0,a_{8}\not=0,a_{12}\not=0$ and $%
a_{2}=a_{3}=a_{5}=a_{6}=a_{7}=a_{9}=a_{10}=a_{11}=a_{13}=0.$
\end{enumerate}

\begin{th}
A conformally flat $Co(QE)_{n}(n>3)$ is a manifold of comprehensive
quasi-constant curvature.
\end{th}

On contracting (\ref{eq-com-cur}), we get 
\begin{eqnarray*}
S(X,Y) &=&b_{1}g(X,Y)+b_{2}\omega ^{1}(X)\omega ^{1}(Y)+b_{3}\omega
^{2}(X)\omega ^{2}(Y) \\
&&+b_{4}\omega ^{3}(X)\omega ^{3}(Y)+b_{5}\omega ^{4}(X)\omega
^{4}(Y)+b_{6}(\omega ^{1}(X)\omega ^{2}(Y)+\omega ^{2}(X)\omega ^{1}(Y)) \\
&&+b_{7}(\omega ^{1}(X)\omega ^{4}(Y)+\omega ^{4}(X)\omega
^{1}(Y))+b_{8}(\omega ^{2}(X)\omega ^{3}(Y)+\omega ^{3}(X)\omega ^{2}(Y)) \\
&&+b_{9}(\omega ^{2}(X)\omega ^{4}(Y)+\omega ^{4}(X)\omega
^{2}(Y))+b_{10}(\omega ^{3}(X)\omega ^{1}(Y)+\omega ^{1}(X)\omega ^{3}(Y)) \\
&&+b_{11}(\omega ^{3}(X)\omega ^{4}(Y)+\omega ^{4}(X)\omega
^{3}(Y))+b_{12}d_{1}(X,Y)+b_{13}d_{2}(X,Y),
\end{eqnarray*}%
where $%
b_{1}=a_{1}(n-1)+a_{4}+a_{5}+a_{6}+a_{7},b_{2}=(n-2)a_{4},b_{3}=(n-2)a_{5},b_{4}=(n-2)a_{6},,b_{5}=(n-2)a_{7},b_{6}=(n-2)a_{8}, 
$ $%
b_{7}=(n-2)a_{9},b_{8}=(n-2)a_{10},b_{9}=(n-2)a_{11},b_{10}=(n-2)a_{12},b_{11}=(n-2)a_{13},b_{12}=(n-2)a_{2},b_{13}=(n-2)a_{3}. 
$ So, we can state that :

\begin{th}
A manifold $M(n>3)$ of comprehensive quasi-constant curvature is a $%
Co(QE)_{n}$.
\end{th}

\section{Existence of $Co(QE)_{n}(n>2)$}

Firstly, we state the the following well known theorem given by Neill \cite%
{Neill}:

\begin{prop}
For a smooth manifold $M,$ the following are equivalent:
\end{prop}

\begin{enumerate}
\item There is a non-vanishing vector field on $M$.

\item Either $M$ is non-compact or compact and has Euler number ${\cal X}%
(M)=0$.
\end{enumerate}

\begin{th}
Let $M$ be a smooth manifold with ${\cal X}(M)=0$. If the Ricci tensor $S$
of the Riemannian manifold is non-vanishing and satisfies 
\begin{eqnarray}
S^{3}(X,Z)S^{3}(Y,W) &=&S(Y,Z)S(X,W)-a_{0}S(X,Z)S(Y,W)  \nonumber \\
&&+a_{1}(S(X,Y)g(Z,W)+S(Z,W)g(X,Y))  \nonumber \\
&&+a_{2}(S^{2}(X,Y)g(Z,W)+S^{2}(Z,W)g(X,Y))  \nonumber \\
&&+a_{3}(S^{3}(X,Y)g(Z,W)+S^{3}(Z,W)g(X,Y))  \nonumber \\
&&+a_{4}(S^{2}(X,Y)S(Z,W)+S^{2}(Z,W)S(X,Y))  \nonumber \\
&&+a_{5}(g(Y,Z)g(X,W)-g(Y,W)g(X,Z))  \nonumber \\
&&+a_{6}(S^{3}(X,Y)S(Z,W)+S^{3}(Z,W)S(X,Y))  \nonumber \\
&&+a_{7}(S^{3}(X,Y)S^{2}(Z,W)+S^{3}(Z,W)S^{2}(X,Y))  \nonumber \\
&&+a_{8}d_{1}(X,W)g(Y,Z)+a_{9}d_{2}(X,W)g(Y,Z)  \nonumber \\
&&+a_{10}S^{2}(X,Z)S^{2}(Y,W),  \label{eq-1-s}
\end{eqnarray}%
where $a_{0},a_{1},\ldots ,a_{10}$ are non-zero scalars and $d_{1},d_{2}$ are
symmetric tensors of type $(0,2)$, then the manifold is $Co(QE)_{n}$.
\end{th}

\noindent {\bf Proof. }Let $M$ be a smooth manifold with ${\cal X}(M)=0$,
then there exist a non-vanishing vector field $U$ on $M$. Let $\omega
^{1}(X)=g(X,U)$ for all vector fields $X$ on $M$. Then taking $Y=Z=U$ in (%
\ref{eq-1-s}) and $\omega ^{2}(X)=\omega ^{1}(QX),\omega ^{3}(X)=\omega
^{1}(Q^{2}X),\omega ^{4}(X)=\omega ^{1}(Q^{3}X)$, we have 
\begin{eqnarray*}
S(U,U)S(X,W) &=&-a_{5}g(U,U)g(X,W)+a_{5}\omega ^{1}(X)\omega ^{1}(W) \\
&&+a_{0}\omega ^{2}(X)\omega ^{2}(W)-a_{10}\omega ^{3}(X)\omega
^{3}(W)+\omega ^{4}(X)\omega ^{4}(W) \\
&&-a_{8}g(U,U)d_{1}(X,W)-a_{9}d_{2}(X,W)g(U,U) \\
&&-a_{1}(\omega ^{2}(X)\omega ^{1}(W)+\omega ^{1}(X)\omega ^{2}(W)) \\
&&-a_{2}(\omega ^{3}(X)\omega ^{1}(W)+\omega ^{1}(X)\omega ^{3}(W)) \\
&&-a_{3}(\omega ^{4}(X)\omega ^{1}(W)+\omega ^{1}(X)\omega ^{4}(W)) \\
&&-a_{4}(\omega ^{3}(X)\omega ^{2}(W)+\omega ^{2}(X)\omega ^{3}(W)) \\
&&-a_{6}(\omega ^{2}(X)\omega ^{4}(W)+\omega ^{4}(X)\omega ^{2}(W)) \\
&&-a_{7}(\omega ^{3}(X)\omega ^{4}(W)+\omega ^{4}(X)\omega ^{3}(W))
\end{eqnarray*}%
which can be written as 
\begin{eqnarray*}
S(X,W) &=&c_{3}g(X,W)+c_{4}\omega ^{1}(X)\omega ^{1}(W)+c_{5}\omega
^{2}(X)\omega ^{2}(W) \\
&&+c_{6}\omega ^{3}(X)\omega ^{3}(W)+c_{7}\omega ^{4}(X)\omega ^{4}(W) \\
&&+c_{8}(\omega ^{2}(X)\omega ^{1}(W)+\omega ^{1}(X)\omega ^{2}(W))) \\
&&+c_{9}(\omega ^{3}(X)\omega ^{1}(W)+\omega ^{1}(X)\omega ^{3}(W)) \\
&&+c_{10}(\omega ^{4}(X)\omega ^{1}(W)+\omega ^{1}(X)\omega ^{4}(W)) \\
&&+c_{11}(\omega ^{3}(X)\omega ^{2}(W)+\omega ^{2}(X)\omega ^{3}(W)) \\
&&+c_{12}(\omega ^{2}(X)\omega ^{4}(W)+\omega ^{4}(X)\omega ^{2}(W)) \\
&&+c_{13}(\omega ^{3}(X)\omega ^{4}(W)+\omega ^{4}(X)\omega ^{3}(W)) \\
&&+c_{14}d_{1}(X,W)+c_{15}d_{2}(X,W),
\end{eqnarray*}%
where 
\begin{eqnarray*}
c_{1} &=&S(U,U),c_{2}=g(U,U),c_{3}=\frac{-a_{5}c_{2}}{c_{1}},c_{4}=\frac{%
a_{5}}{c_{1}}, \\
c_{5} &=&\frac{a_{0}}{c_{1}},c_{6}=-\frac{a_{10}}{c_{1}},c_{7}=\frac{1}{c_{1}%
},c_{8}=-\frac{a_{1}}{c_{1}},c_{9}=-\frac{a_{2}}{c_{1}}, \\
c_{10} &=&-\frac{a_{3}}{c_{1}},c_{11}=-\frac{a_{4}}{c_{1}},c_{12}=-\frac{%
a_{6}}{c_{1}},c_{12}=-\frac{a_{6}}{c_{1}}, \\
c_{13} &=&-\frac{a_{7}}{c_{1}},c_{14}=-\frac{a_{8}c_{2}}{c_{1}},c_{15}=-%
\frac{a_{9}c_{2}}{c_{1}}.
\end{eqnarray*}%
$S(U,U)$ is the Ricci curvature in the direction of the generator $U$ and
Ricci tensor is non-vanishing, so $c_{1}\not=0$. Since $c_{1},\ldots ,c_{15}$
are non-zero scalars. Hence the manifold is $Co(QE)_{n}$.

\section{Sufficient condition for a compact orientable $Co(QE)_{n}$ to be
conformal to a sphere in $\left( n+1\right) $-dimensional Euclidean space}

Now, we state the well known result which is proved by Watanabe \cite[Cor 1]%
{Watanbe}.

In an $n$-dimensional ($n>2$) compact simply connected orientable Riemannian
manifold $M$, we have%
\begin{equation}
\int_{M}S(X,X)dv=\int_{M}\left\vert dX\right\vert ^{2}dv+\frac{n-1}{n}%
\int_{M}\left\vert \delta X\right\vert ^{2}dv  \label{eq-23}
\end{equation}%
for a non-parallel vector field $X$, then the manifold $M$ is conformally
diffeomorphic to a sphere in an $\left( n+1\right) $-dimensional Euclidean
space, where $dv$ is the volume element of $M$ and $dX$ and $\delta X$ are
the curl and divergence of $X$, respectively.

Using this result, we obtain

\begin{th}
Let $M$ be a compact, orientable $Co(QE)_{n}(n\geq 3)$ without boundary and
the generator $W_{1}$ be the gradient of a function. If $W_{1}$ satisfies
the condition 
\[
\int_{M}\left( a+b_{11}\right) dv=\frac{n-1}{n}\int_{M}\left\vert \delta
W_{1}\right\vert ^{2}dv, 
\]%
then the manifold $Co(QE)_{n}$ is conformal to a sphere immersed in
Euclidean space $E^{n+1}$.
\end{th}

\noindent {\bf Proof.} Since $S(W_{1},W_{1})=a+b_{11}$, using (\ref{eq-23}),
we obtain 
\begin{equation}
\int_{M}\left( a+b_{11}\right) dv=\int_{M}\left\vert dW_{1}\right\vert
^{2}dv+\frac{n-1}{n}\int_{M}\left\vert \delta W_{1}\right\vert ^{2}dv.
\label{eq-231}
\end{equation}%
Let $W_{1}={grad}f$, then $dW_{1}=0$. So (\ref{eq-231}) reduces to 
\[
\int_{M}\left( a+b_{11}\right) dv=\frac{n-1}{n}\int_{M}\left\vert \delta
W_{1}\right\vert ^{2}dv. 
\]%
Assume that $W_{1}$ is parallel, then $\nabla W_{1}=0$, that is, $\nabla {%
grad}f=0$ or $\bigtriangleup f=0$, where $\bigtriangleup $ denotes Laplacian
of $f$, and $\nabla $ denotes the covariant differentiation with respect to
the metric of $M$. $\bigtriangleup f=0$ implies that $f=0$ \cite[p. 39]{Yano}%
, therefore $W_{1}=0,$ which contradicts that $W_{1}$ is non-zero. So $W_{1}$
is non-parallel.

Now, by using the Watanabe result, we can say that $Co(QE)_{n}$ is conformal
to a sphere immersed in Euclidean space $E^{n+1}$.

\section{Geometric Properties of $Co(QE)_{n}(n>2)$}

Now, we give some geometrical properties of $Co(QE)_{n}(n>2)$.

\begin{th}
\label{th-1} In a $Co(QE)_{n}(n>2)$, $QW_{1}$ is orthogonal to $W_{1}$ if
and only if $a+b_{11}=0$.
\end{th}

\begin{th}
In a $Co(QE)_{n}(n>2)$, $QW_{1}$ is orthogonal to $W_{2}$ if and only if $%
b_{12}=0$.
\end{th}

\begin{th}
In a $Co(QE)_{n}(n>2)$, $QW_{1}$ is orthogonal to $W_{3}$ if and only if $%
b_{31}=0$.
\end{th}

\begin{th}
In a $Co(QE)_{n}(n>2)$, $QW_{1}$ is orthogonal to $W_{4}$ if and only if $%
b_{41}=0$.
\end{th}

\begin{th}
In a $Co(QE)_{n}(n>2)$, $QW_{2}$ is orthogonal to $W_{2}$ if and only if $%
a+b_{22}+c_{1}d_{1}(W_{2},W_{2})+c_{2}d_{2}(W_{2},W_{2})=0$.
\end{th}

\begin{th}
In a $Co(QE)_{n}(n>2)$, $QW_{3}$ is orthogonal to $W_{3}$ if and only if $%
a+b_{33}+c_{1}d_{1}(W_{3},W_{3})+c_{2}d_{2}(W_{3},W_{3})=0$.
\end{th}

\begin{th}
In a $Co(QE)_{n}(n>2)$, $QW_{4}$ is orthogonal to $W_{4}$ if and only if $%
a+b_{44}+c_{1}d_{1}(W_{4},W_{4})+c_{2}d_{2}(W_{4},W_{4})=0$.
\end{th}

\begin{th}
In a $Co(QE)_{n}(n>2)$, $QW_{2}$ is orthogonal to $W_{3}$ if and only if $%
b_{23}+c_{1}d_{1}(W_{2},W_{3})+c_{2}d_{2}(W_{2},W_{3})=0$.
\end{th}

\begin{th}
In a $Co(QE)_{n}(n>2)$, $QW_{3}$ is orthogonal to $W_{4}$ if and only if $%
b_{34}+c_{1}d_{1}(W_{3},W_{4})+c_{2}d_{2}(W_{3},W_{4})=0$.
\end{th}

\begin{th}
In a $Co(QE)_{n}(n>2)$, $QW_{2}$ is orthogonal to $W_{4}$ if and only if $%
b_{24}+c_{1}d_{1}(W_{2},W_{4})+c_{2}d_{2}(W_{2},W_{4})=0$.
\end{th}

We know that a vector field $X$ on compact orientable Riemannian manifold $M$
without boundary is said to be Killing vector field \cite{Yano} if it
satisfies $\pounds _{X}g=0$. For a Killing vector field $X$, we have the
following result of Yano \cite[p. 43]{Yano} 
\begin{equation}
\int_{M}\left( S(X,X)-\left\vert \nabla X\right\vert ^{2}\right) dv=0,
\label{eq-int}
\end{equation}%
where $dv$ denotes the volume element of $M$.

Let $M$ be a compact orientable $Co(QE)_{n}(n>2)$ without boundary and $X\in
TM$. Let $\alpha$, $\beta$, $\gamma$, $\delta$ be the angle between $W_{1}$ and 
$X$, $W_{2}$ and $X$, $W_{3}$ and $X$, $W_{4}$ and $X$, respectively such
that $\alpha \leq \beta \leq \gamma \leq \delta $. Then $\cos \alpha \geq
\cos \beta \geq \cos \gamma \geq \cos \delta $ and so $g(X,W_{1})\geq
g(X,W_{2})\geq g(X,W_{3})\geq g(X,W_{4}).$ Clearly $g(X,X)\geq \left(
g(X,W_{1})\right) ^{2}.$

From (\ref{eq-co}), we have 
\begin{eqnarray*}
S(X,X) &=&ag(X,X)+b_{ij}\omega ^{i}(X)\omega
^{j}(X)+c_{1}d_{1}(X,X)+c_{2}d_{2}(X,X) \\
&\leq &ag(X,X)+c_{1}d_{1}(X,X)+c_{2}d_{2}(X,X) \\
&&+(b_{11}+b_{22}+b_{33}+b_{44}+2b_{12}+2b_{13}+2b_{23}+2b_{14}+2b_{24}+2b_{34})\left( g(X,W_{1})\right) ^{2}
\\
&\leq
&(a+b_{11}+b_{22}+b_{33}+b_{44}+2b_{12}+2b_{13}+2b_{23}+2b_{14}+2b_{24}+2b_{34})g(X,X),
\end{eqnarray*}%
when $c_{1}d_{1}(X,X)+c_{2}d_{2}(X,X)<0$ and $%
(b_{11}+b_{22}+b_{33}+b_{44}+2b_{12}+2b_{13}+2b_{23}+2b_{14}+2b_{24}+2b_{34})>0 
$. The equation (\ref{eq-int}) gives 
\begin{eqnarray*}
0 &\leq &\int_{M}\left(
(a+b_{11}+b_{22}+b_{33}+b_{44}+2b_{12}+2b_{13}+2b_{23}+2b_{14}+2b_{24}+2b_{34})g(X,X)-\left\vert \nabla X\right\vert ^{2}\right) dv
\\
&=&0.
\end{eqnarray*}%
If $%
(a+b_{11}+b_{22}+b_{33}+b_{44}+2b_{12}+2b_{13}+2b_{23}+2b_{14}+2b_{24}+2b_{34})<0 
$, then we conclude that $g(X,X)=0$ and $\nabla X=0$. Therefore $X=0$. Now,
we can state the following result:

\begin{th}
Let $M$ be a compact orientable $Co(QE)_{n}(n>2)$ without boundary. Then a
Killing vector field other than zero does not exist provided that $%
c_{1}d_{1}(X,X)+c_{2}d_{2}(X,X)$, $%
a+b_{11}+b_{22}+b_{33}+b_{44}+2b_{12}+2b_{13}+2b_{23}+2b_{14}+2b_{24}+2b_{34} 
$ are negative and $%
b_{11}+b_{22}+b_{33}+b_{44}+2b_{12}+2b_{13}+2b_{23}+2b_{14}+2b_{24}+2b_{34}$
is positive.
\end{th}

Now, consider compact orientable conformally flat $Co(QE)_{n}(n>3)$. Let $%
\omega $ be a $p$-form and $F_{p}(\omega ,\omega )$ be a quadratic form \cite%
[p. 70]{Yano} given by 
\begin{equation}
F_{p}(\omega ,\omega )=S_{ji}\omega _{i_{2}\cdots i_{p}}^{j}\omega
^{ii_{2}\cdots i_{p}}+\frac{p-1}{2}R_{kjih}\omega _{i_{3}\cdots
i_{p}}^{kj}\omega ^{ihi_{3}\cdots i_{p}},  \label{eq-p}
\end{equation}%
where $R_{kjih}$ and $S_{ji}$ are the components of the curvature tensor $R$
of type $(0,4)$and the Ricci tensor $S$ of type $(0,2)$ of the $Co(QE)_{n}$.

In virtue of (\ref{eq-co}) and (\ref{eq-con}), we can express (\ref{eq-p})
as follows:

\begin{eqnarray}
F_{p}(\omega ,\omega ) &=&\left( \frac{(1-p)a}{(n-1)(n-2)}+a\right)
\left\vert \omega \right\vert ^{2}+\frac{(p-1)(b_{11}+b_{22}+b_{33}+b_{44})}{%
2(n-1)(n-2)}\left\vert \omega \right\vert ^{2}  \nonumber \\
&&+\frac{\left( n-2p\right) }{(n-2)}b_{lm}\left( W^{l}\cdot \omega \right)
\left( W^{m}\cdot \omega \right)  \nonumber \\
&&+\frac{c_{1}(p-1)}{2(n-2)}\left( 
\begin{array}{c}
2D_{1k}^{h}\omega _{jhi_{3}\cdots i_{p}}\omega ^{kji_{3}\cdots
i_{p}}-D_{1j}^{h}\omega _{khi_{3}\cdots i_{p}}\omega ^{kji_{3}\cdots i_{p}}
\\ 
+D_{1i}^{k}\omega _{hki_{3}\cdots i_{p}}\omega ^{ihi_{3}\cdots i_{p}}%
\end{array}%
\right)  \nonumber \\
&&+\frac{c_{2}(p-1)}{2(n-2)}\left( 
\begin{array}{c}
2D_{2k}^{h}\omega _{jhi_{3}\cdots i_{p}}\omega ^{kji_{3}\cdots
i_{p}}-D_{2j}^{h}\omega _{khi_{3}\cdots i_{p}}\omega ^{kji_{3}\cdots i_{p}}
\\ 
+D_{2i}^{k}\omega _{hki_{3}\cdots i_{p}}\omega ^{ihi_{3}\cdots i_{p}}%
\end{array}%
\right)  \nonumber \\
&&+c_{1}D_{1i}^{h}\omega _{hi_{2}\cdots i_{p}}\omega ^{ii_{2}\cdots
i_{p}}+c_{2}D_{2i}^{h}\omega _{hi_{2}\cdots i_{p}}\omega ^{ii_{2}\cdots
i_{p}},  \label{eq-ppp}
\end{eqnarray}%
where the components of $\omega $ are $\omega _{i_{1}\cdots i_{p}}$ ; $%
W^{l}\cdot \omega $ is a tensor of type $(0,p-1)$ with components $%
W^{li}\omega _{ii_{1}\cdots i_{p-1}}$ and $\left\vert \omega \right\vert
^{2}=\omega _{i_{1}\cdots i_{p}}\omega ^{i_{1}\cdots i_{p}}$.

Let $\omega $ be a Killing $p$-form \cite{Yano}. Then 
\begin{equation}
\int_{Co(QE)_{n}}\left( F_{p}(\omega ,\omega )-\left\vert \omega \right\vert
^{2}\right) dv=0.  \label{eq-pp}
\end{equation}%
In virtue of (\ref{eq-ppp}) and (\ref{eq-pp}), we can say that $\omega =0$
if $c_{1}=0=c_{2}$ and $(p-1)(b_{11}+b_{22}+b_{33}+b_{44}-a)+a(n-1)(n-2)<0$, 
$\left( n-2p\right) <0$ with $l=m$, where $l,m=1,2,3,4$. This leads to the
following result:

\begin{th}
Let $M$ be a compact, orientable conformally flat $Co(QE)_{n}$ $(n>3)$
without boundary. If $c_{1}=0=c_{2}$ and $%
(p-1)(b_{11}+b_{22}+b_{33}+b_{44}-a)+a(n-1)(n-2)<0$, $\left( n-2p\right) <0$
with $l=m$, where $l,m=1,2,3,4$ and $1<p<n$, then there does not exist
non-zero Killing $p$-form.
\end{th}

We give some more results on $Co(QE)_{n}$.

\begin{th}
Let $Co(QE)_{n}$ be a semi-pseudo Ricci symmetric manifold. Then Ricci
tensor of the manifold is not cyclic parallel.
\end{th}

\noindent {\bf Proof. }Let Ricci tensor of semi-pseudo Ricci symmetric $%
Co(QE)_{n}$ be cyclic parallel. Then from (\ref{eq-sp}) and (\ref{eq-cp}),
we have 
\[
\pi (X)S(Y,Z)+\pi (Y)S(X,Z)+\pi (Z)S(X,Y)=0. 
\]%
Using \cite[Lemma 2]{Walker} in above equation, we can say that either all $%
\pi (X)=0$ or all $S(X,Y)=0$. But for $Co(QE)_{n}$, $S(X,Y)\not=0$, so all $%
\pi (X)=0$. But for semi-pseudo Ricci symmetric manifold, $\pi (X)\not=0$.
Therefore our assumption was wrong.

\begin{rem-new}
If the Ricci tensor of semi-pseudo Ricci symmetric $Co(QE)_{n}$ is cyclic
parallel, then it reduces to a Ricci symmetric $Co(QE)_{n}$.
\end{rem-new}

\begin{cor}
A semi-pseudo Ricci symmetric $Co(QE)_{n}$ cannot admit a Codazzi type Ricci
tensor.
\end{cor}

\noindent {\bf Proof.} By using \cite[Th 5]{Tarafdar}, we get the result.

\begin{th}
If the generator $W_{1}$ of a $Co(QE)_{n}$ is a concurrent vector field,
then $QW_{1}$ is orthogonal to $W_{1}$.
\end{th}

\noindent {\bf Proof.} Let the generator $W_{1}$ of a $Co(QE)_{n}$ be a concurrent
vector field. Then it is easy to verify that $R(X,Y)W_{1}=0$, so $%
S(Y,W_{1})=0$. Using (\ref{eq-co}), we have 
\[
0=S(X,W_{1})=ag(X,W_{1})+b_{11}\omega ^{1}(X)+b_{21}\omega
^{2}(X)+b_{31}\omega ^{3}(X)+b_{41}\omega ^{4}(X). 
\]%
Taking $X=W_{1}$ in above, we get $a+b_{11}=0$. By using Theorem \ref{th-1},
we get the result.

\begin{th}
Let the generator $W_{1}$ of a $Co(QE)_{n}$ be a concircular vector field
and the associated scalars be constants. Then the associated $1$-forms $%
\omega ^{1},\omega ^{2},\omega ^{3},\omega ^{4}$ are closed provided $%
a+b_{11}\not=0$, $b_{21}$,$b_{31},b_{41}\not=0$.
\end{th}

\noindent {\bf Proof.} Let the generator $W_{1}$ of a $Co(QE)_{n}$ be a
concircular vector field. Then $R(X,Y)W_{1}=\left( X\rho \right) Y-\left(
Y\rho \right) X$, so $S(Y,W_{1})=(1-n)\left( Y\rho \right) $.

Using (\ref{eq-co}), we have 
\[
(1-n)\left( XY\rho \right) =\left( a+b_{11}\right) \nabla _{X}\omega
^{1}(Y)+b_{21}\nabla _{X}\omega ^{2}(Y)+b_{31}\nabla _{X}\omega
^{3}(Y)+b_{41}\nabla _{X}\omega ^{4}(Y), 
\]%
\[
(1-n)\left( YX\rho \right) =\left( a+b_{11}\right) \nabla _{Y}\omega
^{1}(X)+b_{21}\nabla _{Y}\omega ^{2}(X)+b_{31}\nabla _{Y}\omega
^{3}(X)+b_{41}\nabla _{Y}\omega ^{4}(X), 
\]%
\[
(1-n)\left( [X,Y]\rho \right) =\left( a+b_{11}\right) \omega
^{1}([X,Y])+b_{21}\omega ^{2}([X,Y])+b_{31}\omega ^{3}([X,Y])+b_{41}\omega
^{4}([X,Y]), 
\]%
by using these three equations, we have 
\begin{eqnarray*}
0 &=&\left( a+b_{11}\right) \left( \left( \nabla _{X}\omega ^{1}\right)
(Y)-\left( \nabla _{Y}\omega ^{1}\right) (X)\right) +b_{21}\left( \left(
\nabla _{X}\omega ^{2}\right) (Y)-\left( \nabla _{Y}\omega ^{2}\right)
(X)\right) \\
&&+b_{31}\left( \left( \nabla _{X}\omega ^{3}\right) (Y)-\left( \nabla
_{Y}\omega ^{3}\right) (X)\right) +b_{41}\left( \left( \nabla _{X}\omega
^{4}\right) (Y)-\left( \nabla _{Y}\omega ^{4}\right) (X)\right) .
\end{eqnarray*}

\begin{th}
Let the associated $1$-form $\omega ^{1}$ in a $Co(QE)_{n}$ be closed. Then
the integral curves of the vector field $W_{1}$ are geodesic.
\end{th}

\noindent {\bf Proof.} The proof is similar to \cite[Th 2.3]{Mallick}.

\begin{cor}
Let the generator $W_{1}$ of a $Co(QE)_{n}$ be a concircular vector field.
Then the integral curves of the vector field $W_{1}$ are geodesic, provided $%
a+b_{11}\not=0$.
\end{cor}

\begin{th}
Let the generators of a $Co(QE)_{n}$ be Killing vector fields and the
associated scalars be constants. Then the Ricci tensor of the manifold is
cyclic parallel if and only if the structure tensors are cyclic parallel.
\end{th}

\noindent {\bf Proof.} Let the generators of a $Co(QE)_{n}$ be Killing vector
fields, that is, $\left( \nabla _{X}\omega ^{i}\right) (Y)+\left( \nabla
_{Y}\omega ^{i}\right) (X)=0$, $i=1,2,3,4$. Then it is easy to calculate
that 
\begin{eqnarray*}
&&\left( \nabla _{X}S\right) (Y,Z)+\left( \nabla _{Y}S\right) (Z,X)+\left(
\nabla _{Z}S\right) (X,Y) \\
&=&c_{1}\left( \left( \nabla _{X}d_{1}\right) (Y,Z)+\left( \nabla
_{Y}d_{1}\right) (Z,X)+\left( \nabla _{Z}d_{1}\right) (X,Y)\right) \\
&&+c_{2}\left( \left( \nabla _{X}d_{2}\right) (Y,Z)+\left( \nabla
_{Y}d_{2}\right) (Z,X)+\left( \nabla _{Z}d_{2}\right) (X,Y)\right) .
\end{eqnarray*}

\begin{th}
If the generators of a $Co(QE)_{n}$ corresponding to the associated $1$%
-forms are recurrent with the same vector of recurrence and the associated
scalars are constants, then the manifold is a generalized Ricci recurrent
manifold provided that $c_{1}=0=c_{2}.$
\end{th}

\noindent {\bf Proof.} Let the generators $W_{i}$ of a $Co(QE)_{n}$ corresponding
to the associated $1$-forms $\omega ^{i}$, $i=1,2,3,4$., respectively, be
recurrent with the same vector of recurrence, that is, $\left( \nabla
_{X}\omega ^{i}\right) (Y)=\alpha (X)\omega ^{i}(Y),$ $i=1,2,3,4,$ where
where $\alpha (X)$ is a nonzero $1$-form. Using (\ref{eq-S}), we obtain 
\begin{eqnarray*}
\left( \nabla _{W}S\right) (X,Y) &=&2\alpha (W)\left(
S(X,Y)-ag(X,Y)-c_{1}d_{1}(X,Y)+c_{2}d_{2}(X,Y)\right) \\
&&-c_{1}\left( \nabla _{W}d_{1}\right) (X,Y)-c_{2}\left( \nabla
_{W}d_{1}\right) (X,Y).
\end{eqnarray*}

\section{Sectional curvatures at a point of a conformally flat $Co(QE)_{n}$}

Let $\{W_{1},W_{2},W_{3},W_{4}\}^{\bot }$ denote the $(n-4)$-dimensional
distribution in a conformally flat $Co(QE)_{n}$ $(n>3)$ orthogonal to $%
\{W_{1},W_{2},W_{3},W_{4}\}$. Let $X,Y\in \{W_{1},W_{2},W_{3},W_{4}\}^{\bot
} $. By using (\ref{eq-con}), we have 
\begin{eqnarray*}
R(X,Y,Y,X) &=&\frac{1}{n-2}\left(
S(X,X)g(Y,Y)+S(Y,Y)g(X,X)-2S(X,Y)g(X,Y)\right) \\
&&-\frac{r}{(n-1)(n-2)}\left( g(X,X)g(Y,Y)-g(Y,X)g(X,Y)\right) .
\end{eqnarray*}%
By using (\ref{eq-co}) and taking $c_{1}=0=c_{2}$, we get the sectional
curvature of the plane determined by two vectors $X,Y\in
\{W_{1},W_{2},W_{3},W_{4}\}^{\bot }$ is 
\[
K(X,Y)=\frac{R(X,Y,Y,X)}{g(X,X)g(Y,Y)-\left( g(X,Y)\right) ^{2}}=\frac{%
a(n-2)-b_{11}-b_{22}-b_{33}-b_{44}}{(n-1)(n-2)}. 
\]%
Assume that $c_{1}=0=c_{2}.$ By similar process, it is easy to calculate
that the sectional curvature of the plane determined by two vectors $X,\in
\{W_{1},W_{2},W_{3},W_{4}\}^{\bot }$ and $W_{1}$; $X,\in
\{W_{1},W_{2},W_{3},W_{4}\}^{\bot }$ and $W_{2}$; $X,\in
\{W_{1},W_{2},W_{3},W_{4}\}^{\bot }$ and $W_{3}$; $X,\in
\{W_{1},W_{2},W_{3},W_{4}\}^{\bot }$ and $W_{4}$ are%
\[
K(X,W_{1})=\frac{R(X,W_{1},W_{1},X)}{g(X,X)}=\frac{\left( a+b_{11}\right)
(n-2)-b_{22}-b_{33}-b_{44}}{(n-1)(n-2)}, 
\]%
\[
K(X,W_{2})=\frac{R(X,W_{2},W_{2},X)}{g(X,X)}=\frac{\left( a+b_{22}\right)
(n-2)-b_{11}-b_{33}-b_{44}}{(n-1)(n-2)}, 
\]%
\[
K(X,W_{1})=\frac{R(X,W_{3},W_{3},X)}{g(X,X)}=\frac{\left( a+b_{33}\right)
(n-2)-b_{11}-b_{22}-b_{44}}{(n-1)(n-2)}, 
\]%
\[
K(X,W_{1})=\frac{R(X,W_{4},W_{4},X)}{g(X,X)}=\frac{\left( a+b_{44}\right)
(n-2)-b_{11}-b_{22}-b_{33}}{(n-1)(n-2)}, 
\]%
respectively.

\section{General two viscous fluid $Co(QE)_{4}$ spacetime}

A viscous fluid spacetime $(M,g)$ is a connected $4$-dimensional
semi-Riemannian manifold with Lorentzian metric $g$ of signature $(-,+,+,+)$%
. The $(0,2)$-type energy momentum tensor $T$ in a general two viscous fluid
spacetime \cite{Coley} is of the form 
\begin{eqnarray}
T(X,Y) &=&p_{r}g(X,Y)+(\sigma _{r}+p_{r})\omega ^{r}(X)\omega
^{r}(Y)-\varsigma _{r}e_{r}(X,Y)+q^{r}(X)\omega ^{r}(Y)\nonumber \\&&+q^{r}(Y)\omega
^{r}(X)  
+(\sigma _{m}+p_{m})\omega ^{m}(X)\omega ^{m}(Y)+p_{m}g(X,Y)-\varsigma
_{m}e_{m}(X,Y)\nonumber \\&&+q^{m}(X)\omega ^{m}(Y)+q^{m}(Y)\omega ^{m}(X)  \label{eq-p1}
\end{eqnarray}%
together with $g(X,W_{r})=\omega ^{r}(X),g(X,W_{m})=\omega
^{m}(X),g(X,Q_{r})=q^{r}(X),g(X,Q_{m})=q^{m}(X)$ such that $\omega
^{m}(W_{m})=-1$, $\omega ^{r}(W_{r})=-1,q^{r}(Q_{r})=1,q^{m}(Q_{m})=1,\omega
^{m}(W_{r})=0,q^{r}(Q_{m})=0,q^{r}(W_{m})=0,\omega ^{r}(Q_{m})=0$, where $%
\sigma _{r}$, $\sigma _{m}$ are the energy density, $p_{r},p_{m}$ the
isotropic pressure, $\varsigma _{r},\varsigma _{m}$ the shear viscosity
coefficient, $e_{r}$, $e_{m}$ the shear tensor, $\omega _{r},\omega _{m}$
the velocities of the radiation and matter fields and $q_{r},q_{m}$ the heat
conduction vector field in the two viscous fluid. Then in the general
relativistic spacetime whose matter content is viscous fluid obeying the
Einstein's field equation, the Ricci tensor satisfies the following equation 
\begin{equation}
S(X,Y)-\frac{r}{2}g(X,Y)+\Lambda g(X,Y)=\kappa T(X,Y),  \label{eq-p2}
\end{equation}%
where $r$ is the scalar curvature, $\kappa $ is the cosmological constant, $%
\Lambda $ is the gravitational constant. By using (\ref{eq-p1}), (\ref{eq-p2}%
) reduces to 
\begin{eqnarray*}
S(X,Y)-\frac{r}{2}g(X,Y)+\Lambda g(X,Y) &=&\kappa \left( 
\begin{array}{c}
p_{r}g(X,Y)+(\sigma _{r}+p_{r})\omega ^{r}(X)\omega ^{r}(Y) \\ 
-\varsigma _{r}e_{r}(X,Y)+q^{r}(X)\omega ^{r}(Y)+q^{r}(Y)\omega ^{r}(X)%
\end{array}%
\right) \\
&&+\kappa \left( 
\begin{array}{c}
(\sigma _{m}+p_{m})\omega ^{m}(X)\omega ^{m}(Y)+p_{m}g(X,Y) \\ 
-\varsigma _{m}e_{m}(X,Y)+q^{m}(X)\omega ^{m}(Y)+q^{m}(Y)\omega ^{m}(X)%
\end{array}%
\right) ,
\end{eqnarray*}%
which gives 
\begin{eqnarray}
S(X,Y) &=&\left( \kappa p_{r}+\kappa p_{m}-\Lambda +\frac{r}{2}\right)
g(X,Y)+\kappa (\sigma _{r}+p_{r})\omega ^{r}(X)\omega ^{r}(Y)  \nonumber \\
&&+\kappa (\sigma _{m}+p_{m})\omega ^{m}(X)\omega ^{m}(Y)+\kappa
q^{r}(X)\omega ^{r}(Y)  \nonumber \\
&&+\kappa q^{r}(Y)\omega ^{r}(X)+\kappa q^{m}(X)\omega ^{m}(Y)+\kappa
q^{m}(Y)\omega ^{m}(X)  \nonumber \\
&&-\kappa \varsigma _{r}e_{r}(X,Y)-\kappa \varsigma _{m}e_{m}(X,Y).
\label{eq-TTT}
\end{eqnarray}

Hence, we can state the following:

\begin{th}
A general two viscous fluid spacetime admitting heat flux and obeying
Einstein's field equation with a cosmological constant is a $Co(QE)_{4}$
spacetime.
\end{th}

\begin{th}
Let $Co(QE)_{4}$ be general two viscous fluid spacetime which is admitting
heat flux and satisfies Einstein's field equation with cosmological
constant, then energy densities of the fluid cannot be a constant.
\end{th}

\noindent {\bf Proof.} Assuming $\omega ^{1}$ as $\omega ^{r}$, $\omega ^{2}$
as $\omega ^{m}$, $\omega ^{3}$ as $q^{r}$, $\omega ^{4}$ as $q^{m}$, using (%
\ref{eq-co}) in (\ref{eq-TTT}) and taking $X=Y=W_{r}$, we get 
\[
\sigma _{r}=\frac{a+b_{11}-b_{22}-b_{33}-b_{44}+2\Lambda }{2\kappa }-\left(
2p_{r}+p_{m}\right) . 
\]%
Using (\ref{eq-co}) in (\ref{eq-TTT}) and taking $X=Y=W_{m}$, we get 
\[
\sigma _{m}=\frac{a+b_{22}-b_{11}-b_{33}-b_{44}+2\Lambda }{2\kappa }-\left(
2p_{m}+p_{r}\right) . 
\]

\section{$Co(QE)_{n}$ spacetime with vanishing space-matter \\tensor}

Let $(M,g)$ be a smooth manifold. Petrov \cite{Petrov} introduced a $(0,4)$%
-type tensor $P$, which is known as the space-matter tensor of the manifold
and defined as%
\begin{equation}
P=R+\frac{\kappa }{2}g\wedge T-\sigma G,  \label{eq-P1}
\end{equation}%
where $R$ is the curvature tensor of type $(0,4)$, $T$ is the energy
momentum tensor of type $(0,2)$, $\kappa $ is the gravitational constant, $%
\sigma $ is the energy density, $G$ is a tensor of type $(0,4)$ given by 
\[
G(X,Y,Z,W)=g(Y,Z)g(X,W)-g(X,Z)g(Y,W)
\]%
and Kulkarni--Nomizu product $\alpha \wedge \beta $ of two $(0,2)$ tensors $%
\alpha $ and $\beta $ is defined by 
\begin{eqnarray*}
(\alpha \wedge \beta )(X,Y,Z,W) &=&\alpha (Y,Z)\beta (X,W)+\alpha (X,W)\beta
(Y,Z) \\
&&-\alpha (X,Z)\beta (Y,W)-\alpha (Y,W)\beta (X,Z).
\end{eqnarray*}%
Let $P=0$ in (\ref{eq-P1}) and using (\ref{eq-p2}), we get 
\[
R=-\frac{1}{2}g\wedge \left( S-\frac{r}{2}g+\Lambda g\right) +\sigma G.
\]%
Then using (\ref{eq-p2}) and (\ref{eq-P1}), we have 
\begin{eqnarray}
\left( {div}P\right) (X,Y,Z) &=&\left( {div}R\right) (X,Y,Z)+\frac{1}{2}%
\left( \left( \nabla _{X}S\right) (Y,Z)-\left( \nabla _{Y}S\right)
(X,Z)\right)   \nonumber \\
&&-g(Y,Z)\left( d\sigma (X)+\frac{1}{4}dr(X)\right) +g(X,Z)\left( d\sigma
(Y)+\frac{1}{4}dr(Y)\right)   \nonumber \\
&=&\frac{3}{2}\left( \left( \nabla _{X}S\right) (Y,Z)-\left( \nabla
_{Y}S\right) (X,Z)\right)   \nonumber \\
&&-g(Y,Z)\left( d\sigma (X)+\frac{1}{4}dr(X)\right) +g(X,Z)\left( d\sigma
(Y)+\frac{1}{4}dr(Y)\right) .  \label{eq-P}
\end{eqnarray}

\begin{th}
In a $Co(QE)_{4}$ spacetime satisfying Einstein's field equation with
divergence free space-matter tensor the energy density is constant.
\end{th}

\noindent {\bf Proof.} Consider (\ref{eq-P}) with ${div}P=0$. On contracting
(\ref{eq-P}) over $Y$ and $Z$, we have $d\sigma (X)=0$.

\begin{th}
Let the associated scalars and the energy density in a $Co(QE)_{4}$
spacetime satisfying Einstein's field equation be constants and generators
of $Co(QE)_{4}$ be parallel. Then the space-matter tensor will be divergence
free.
\end{th}

\noindent {\bf Proof.} Using (\ref{eq-co}) in (\ref{eq-P}) and taking the
associated scalars and the energy density of $Co(QE)_{4}$ as constants with
parallel generators. Then we get ${div}P=0$.

\section{Example of $Co(QE)_{4}$}

Consider the G\"{o}del metric \cite{Godel}, which is defined in a $4$%
-dimensional manifold $M$ 
\[
ds^{2}=k^{2}\left( (dt+e^{x}dy)^{2}-dx^{2}-\frac{e^{2x}}{2}%
dy^{2}-dz^{2}\right) ,
\]%
where $k$ is real. The non-vanishing components of metric tensor are 
\[
g_{11}=k^{2},g_{22}=-k^{2},g_{33}=\frac{k^{2}e^{2x}}{2}%
,g_{44}=-k^{2},g_{13}=g_{31}=k^{2}e^{x}.
\]%
The non-vanishing components of inverse of metric tensor are 
\[
g^{11}=g^{22}=-\frac{1}{k^{2}},g^{33}=-\frac{2e^{-2x}}{k^{2}},g^{44}=-\frac{1%
}{k^{2}},g^{13}=g^{31}=\frac{2e^{-x}}{k^{2}}.
\]%
The non-zero Christoffel symbols are 
\[
\Gamma _{12}^{1}=1,\Gamma _{23}^{1}=\Gamma _{13}^{2}=\frac{e^{x}}{2},\Gamma
_{33}^{2}=\frac{e^{2x}}{2},\Gamma _{12}^{3}=-e^{-x}.
\]%
The non-vanishing components of Ricci tensor are 
\[
S_{11}=1,S_{13}=e^{x},S_{33}=e^{2x}
\]%
and the scalar curvature is 
\[
r=-\frac{1}{k^{2}}.
\]%
Now, we choose the functions as follows: 
\[
a=-\frac{1}{k^{2}},\quad b_{11}=-\frac{1}{k^{2}},\quad
b_{12}=b_{13}=b_{24}=b_{14}=0,\quad b_{34}=\frac{2\sqrt{2}}{k^{2}},
\]%
\[
b_{23}=\frac{\sqrt{2}}{k^{2}},\quad b_{22}=-\frac{3}{2k^{2}},\quad b_{33}=%
\frac{3}{k^{2}},\quad b_{44}=\frac{5}{2k^{2}},\quad
c_{1}=e^{-2x},c_{2}=-e^{2x}.
\]%
We take the $1$-forms as follows:%
\[
\omega _{i}^{1}(p)=\left\{ 
\begin{array}{cc}
k & i=4 \\ 
0 & {\rm otherwise}%
\end{array}%
\right. ,
\]%
\[
\omega _{i}^{2}(p)=\left\{ 
\begin{array}{cc}
k & i=2 \\ 
0 & {\rm otherwise}%
\end{array}%
\right. ,
\]%
\[
\omega _{i}^{3}(p)=\left\{ 
\begin{array}{cc}
\dfrac{ke^{x}}{\sqrt{2}} & i=3 \\ 
0 & {\rm otherwise}%
\end{array}%
\right. ,
\]%
\[
\omega _{i}^{4}(p)=\left\{ 
\begin{array}{cc}
k & i=1 \\ 
0 & {\rm otherwise}%
\end{array}%
\right. 
\]%
at any ponit $p\in M$. We choose the associated tensors as follows: 
\[
d_{1}(p)=\left[ 
\begin{array}{cccc}
-\dfrac{e^{2x}}{4} & e^{2x} & 0 & 0 \\ 
e^{2x} & \dfrac{3e^{2x}}{4} & 0 & 0 \\ 
0 & 0 & -\dfrac{e^{2x}}{2} & 0 \\ 
0 & 0 & 0 & 0%
\end{array}%
\right] ,
\]%
\[
d_{2}(p)=\left[ 
\begin{array}{cccc}
\dfrac{e^{-2x}}{4} & e^{-2x} & 0 & 0 \\ 
e^{-2x} & \dfrac{e^{-2x}}{4} & e^{-x} & 0 \\ 
0 & e^{-x} & -\dfrac{e^{-2x}}{2} & 0 \\ 
0 & 0 & 0 & 0%
\end{array}%
\right] 
\]%
at any ponit $p\in M$. Clearly, we can check that the trace of $(0,2)$%
-tensors $d_{1}$ and $d_{2}$ are zero. It is easy to verify that $1$-forms $%
\omega _{i}^{1},\omega _{i}^{2},\omega _{i}^{3},\omega _{i}^{4}$ are unit
and orthogonal. We can see that $d_{1}^{ij}X_{i}W_{j}^{1}=0$ and $%
d_{2}^{ij}X_{i}W_{j}^{1}=0$ for any vector field $X$ on $M$. Therefore, we
can say that manifold with G\"{o}del metric is a $Co(QE)_{4}$ spacetime.

\noindent {\bf Acknowledgement:} The first author would like to thank Prof. Mukut
Mani Tripathi for useful suggestions. 

\noindent {\bf Declaration: }

\noindent{\bf Data Availability: }No data were used to support this study.

\noindent{\bf Conflicts of Interest: }The authors declare that they have no conflicts
of interest.

Department of Mathematics \& Statistics

School of Mathematical \& Physical Sciences

Dr. Harisingh Gour University

Sagar- 470 003, MP, INDIA

Email: {\em pgupta@dhsgsu.edu.in}
\bigskip

and

\bigskip

Department of Mathematics

Indian Institute of Science Education \& Research

Bhopal- 462 066, MP, INDIA

Email: {\em sanjayks@iiserb.ac.in}

\end{document}